\documentclass{amsart}
\usepackage{amssymb,amsmath}
\usepackage[mathscr]{euscript}
\input{epsf}

\usepackage{bbding}
\usepackage{textcomp}

\newcounter{sec}

\newcounter{punct}[sec]

\def\punct{\refstepcounter{punct}{\arabic{sec}.%
\arabic{punct}.  }}

\newtheorem{theorem}{Theorem}[sec]

\newtheorem{definition}[theorem]{Definition}

\newtheorem{conjecture}[theorem]{Conjecture}

\def\COUNTERS{\addtocounter{sec}{1}
              \setcounter{punct}{0}
          \setcounter{equation}{0}
          \setcounter{theorem}{0}
          }
          
          \def\sm{\smallskip}

\begin{document}

\newcommand{\supp}{\mathop {\mathrm {supp}}\nolimits}
\newcommand{\rk}{\mathop {\mathrm {rk}}\nolimits}
\newcommand{\Aut}{\mathop {\mathrm {Aut}}\nolimits}
\newcommand{\Ob}{\mathop {\mathrm {Ob}}\nolimits}
\newcommand{\Out}{\mathop {\mathrm {Out}}\nolimits}
\renewcommand{\Re}{\mathop {\mathrm {Re}}\nolimits}
\newcommand{\Inn}{\mathop {\mathrm {Inn}}\nolimits}
\newcommand{\Char}{\mathop {\mathrm {Char}}\nolimits}
\newcommand{\ch}{\cosh}
\newcommand{\sh}{\sinh}
\newcommand{\Sp}{\mathop {\mathrm {Sp}}\nolimits}
\newcommand{\SOS}{\mathop {\mathrm {SO^*}}\nolimits}
\newcommand{\Ams}{\mathop {\mathrm {Ams}}\nolimits}
\newcommand{\Gms}{\mathop {\mathrm {Gms}}\nolimits}

\def\0{\mathbf 0}

\def\ov{\overline}
\def\wh{\widehat}
\def\wt{\widetilde}
\def\pol{\twoheadrightarrow}

\renewcommand{\rk}{\mathop {\mathrm {rk}}\nolimits}
\renewcommand{\Aut}{\mathop {\mathrm {Aut}}\nolimits}
\renewcommand{\Re}{\mathop {\mathrm {Re}}\nolimits}
\renewcommand{\Im}{\mathop {\mathrm {Im}}\nolimits}
\newcommand{\sgn}{\mathop {\mathrm {sgn}}\nolimits}

\def\bfa{\mathbf a}
\def\bfb{\mathbf b}
\def\bfc{\mathbf c}
\def\bfd{\mathbf d}
\def\bfe{\mathbf e}
\def\bff{\mathbf f}
\def\bfg{\mathbf g}
\def\bfh{\mathbf h}
\def\bfi{\mathbf i}
\def\bfj{\mathbf j}
\def\bfk{\mathbf k}
\def\bfl{\mathbf l}
\def\bfm{\mathbf m}
\def\bfn{\mathbf n}
\def\bfo{\mathbf o}
\def\bfp{\mathbf p}
\def\bfq{\mathbf q}
\def\bfr{\mathbf r}
\def\bfs{\mathbf s}
\def\bft{\mathbf t}
\def\bfu{\mathbf u}
\def\bfv{\mathbf v}
\def\bfw{\mathbf w}
\def\bfx{\mathbf x}
\def\bfy{\mathbf y}
\def\bfz{\mathbf z}

\def\bfA{\mathbf A}
\def\bfB{\mathbf B}
\def\bfC{\mathbf C}
\def\bfD{\mathbf D}
\def\bfE{\mathbf E}
\def\bfF{\mathbf F}
\def\bfG{\mathbf G}
\def\bfH{\mathbf H}
\def\bfI{\mathbf I}
\def\bfJ{\mathbf J}
\def\bfK{\mathbf K}
\def\bfL{\mathbf L}
\def\bfM{\mathbf M}
\def\bfN{\mathbf N}
\def\bfO{\mathbf O}
\def\bfP{\mathbf P}
\def\bfQ{\mathbf Q}
\def\bfR{\mathbf R}
\def\bfS{\mathbf S}
\def\bfT{\mathbf T}
\def\bfU{\mathbf U}
\def\bfV{\mathbf V}
\def\bfW{\mathbf W}
\def\bfX{\mathbf X}
\def\bfY{\mathbf Y}
\def\bfZ{\mathbf Z}

\def\frD{\mathfrak D}
\def\frQ{\mathfrak Q}
\def\frS{\mathfrak S}
\def\frT{\mathfrak T}
\def\frL{\mathfrak L}
\def\frM{\mathfrak M}
\def\frG{\mathfrak G}
\def\frb{\mathfrak b}
\def\frg{\mathfrak g}
\def\frh{\mathfrak h}
\def\frf{\mathfrak f}
\def\frk{\mathfrak k}
\def\frl{\mathfrak l}
\def\frm{\mathfrak m}
\def\frn{\mathfrak n}
\def\fro{\mathfrak o}
\def\frp{\mathfrak p}
\def\frq{\mathfrak q}
\def\frr{\mathfrak r}
\def\frs{\mathfrak s}
\def\frt{\mathfrak t}
\def\fru{\mathfrak u}
\def\frv{\mathfrak v}
\def\frw{\mathfrak w}
\def\frx{\mathfrak x}
\def\fry{\mathfrak y}
\def\frz{\mathfrak z}

\def\bfw{\mathbf w}

\def\R {{\mathbb R }}
 \def\C {{\mathbb C }}
  \def\Z{{\mathbb Z}}
  \def\H{{\mathbb H}}
\def\K{{\mathbb K}}
\def\N{{\mathbb N}}
\def\Q{{\mathbb Q}}
\def\A{{\mathbb A}}

\def\T{\mathbb T}
\def\P{\mathbb P}
\def\SS{\mathbb S}

\def\G{\mathbb G}

\def\cD{\EuScript D}
\def\cL{\EuScript L}
\def\cK{\EuScript K}
\def\cM{\EuScript M}
\def\cN{\EuScript N}
\def\cP{\EuScript P}
\def\cT{\EuScript T}
\def\cQ{\EuScript Q}
\def\cR{\EuScript R}
\def\cW{\EuScript W}
\def\cY{\EuScript Y}
\def\cF{\EuScript F}
\def\cG{\EuScript G}
\def\cZ{\EuScript Z}
\def\cI{\EuScript I}
\def\cB{\EuScript B}
\def\cA{\EuScript A}
\def\cE{\EuScript E}
\def\cC{\EuScript C}
\def\cS{\EuScript S}

\def\bbA{\mathbb A}
\def\bbB{\mathbb B}
\def\bbD{\mathbb D}
\def\bbE{\mathbb E}
\def\bbF{\mathbb F}
\def\bbG{\mathbb G}
\def\bbI{\mathbb I}
\def\bbJ{\mathbb J}
\def\bbL{\mathbb L}
\def\bbM{\mathbb M}
\def\bbN{\mathbb N}
\def\bbO{\mathbb O}
\def\bbP{\mathbb P}
\def\bbQ{\mathbb Q}
\def\bbS{\mathbb S}
\def\bbT{\mathbb T}
\def\bbU{\mathbb U}
\def\bbV{\mathbb V}
\def\bbW{\mathbb W}
\def\bbX{\mathbb X}
\def\bbY{\mathbb Y}

\def\kappa{\varkappa}
\def\epsilon{\varepsilon}
\def\phi{\varphi}
\def\le{\leqslant}
\def\ge{\geqslant}

\def\B{\mathrm B}

\def\la{\langle}
\def\ra{\rangle}
\def\tri{\triangleright}

\def\lambdA{{\boldsymbol{\lambda}}}
\def\alphA{{\boldsymbol{\alpha}}}
\def\betA{{\boldsymbol{\beta}}}
\def\mU{{\boldsymbol{\mu}}}

\def\const{\mathrm{const}}
\def\rem{\mathrm{rem}}
\def\even{\mathrm{even}}
\def\SO{\mathrm{SO}}
\def\SL{\mathrm{SL}}
\def\PSL{\mathrm{PSL}}
\def\cont{\mathrm{cont}}
\def\O{\mathrm{O}}

\def\U{\operatorname{U}}
\def\GL{\operatorname{GL}}
\def\Mat{\operatorname{Mat}}
\def\End{\operatorname{End}}
\def\Mor{\operatorname{Mor}}
\def\Aut{\operatorname{Aut}}
\def\inv{\operatorname{inv}}
\def\red{\operatorname{red}}
\def\Ind{\operatorname{Ind}}
\def\dom{\operatorname{dom}}
\def\im{\operatorname{im}}
\def\md{\operatorname{mod\,}}
\def\indef{\operatorname{indef}}
\def\Gr{\operatorname{Gr}}
\def\Pol{\operatorname{Pol}}
\def\Cut{\operatorname{Cut}}
\def\Add{\operatorname{Add}}
\def\ord{\operatorname{ord}}
\def\Replace{\operatorname{Replace}}
\def\Tr{\operatorname{Tr}}

\def\arr{\rightrightarrows}
\def\bs{\backslash}

\def\cH{\EuScript{H}}
\def\cO{\EuScript{O}}
\def\cQ{\EuScript{Q}}
\def\cL{\EuScript{L}}
\def\cX{\EuScript{X}}

\def\Di{\Diamond}
\def\di{\diamond}

\def\fin{\mathrm{fin}}
\def\ThetA{\boldsymbol {\Theta}}

\def\krest{{\text{\tiny{\rm\texttimes}}}}
\def\vph{{\vphantom{\bigr|}}}

 \def\dow{{\left.\vphantom{\bigl|}\right\downarrow}}
 
 \def\cir{{\circledcirc}}

\begin{center}

\Large \bf

Non-singular actions of infinite-dimensional groups and polymorphisms

\medskip

\sc Yury A.Neretin%
	\footnote{Supported by the grant FWF  P31591.}

\end{center}

\bigskip

{\small Let $Z$ be a probabilistic measure space with a measure $\zeta$, $\mathbb{R}^\times$ be the multiplicative group of positive reals, let $t$ be the coordinate on $\mathbb{R}^\times$.
A polymorphism of $Z$ is a measure $\pi$ on $Z\times Z\times \mathbb{R}^\times$ such that
for any measurable $A$, $B\subset Z$ we have $\pi(A\times Z\times \mathbb{R}^\times)=\zeta(A)$ and the integral $\int t\,d\pi(z,u,t)$ 
over 
$Z\times B\times \mathbb{R}^\times$ is $\zeta(B)$. The set  of all polymorphisms has a natural semigroup structure,
  the group of all nonsingular transformations is dense in this semigroup. 
We discuss a problem of closure in polymorphisms  for certain types of infinite dimensional ('large') groups 
and show that
 a non-singular action of an infinite-dimensional group
generates a representation of its train (category of double cosets) by polymorphisms.

}

\bigskip

\section{Polymorphisms}

\COUNTERS

{\bf \punct Introduction.}
Let $Z=(Z,\zeta)$ be a Lebesgue probabilistic measure space with a non-atomic measure.  
 It is well-known that the  group of all measure preserving transformations of $Z$ has a natural compactification, which is a
metrizable semigroup with a separately continuous multiplication
(this topic arises to E.~Hopf, \cite{Hopf}, 1954).
Its points (we use the term '{\it polymorphism}' proposed in \cite{Ver1}) are Borel measures $\pi$ on $Z\times Z$, whose pushforwards
under the projections of $Z\times Z$ to the first and second factors coincide with $\zeta$, for details, see, e.g. 
\cite{Kre}, \cite{Ver1}, \cite{Ver2}, \cite{Ner-book}, Sect. VIII.4.  Such objects play
 different roles in ergodic theory. They are  a kind of multivalued transformations
of measure spaces and have the same rights as usual measure preserving transformations, see, e.g., \cite{Kre}, \cite{Ver1}. On the other hand, they can be intertwiners between dynamical systems (joinings), see, e.g., 
\cite{Rud}, \cite{King}, \cite{Gla}, \cite{Ryzh}. 

Closures of groups in polymorphisms is an interesting and non-trivial topic,  
starting from the group $\Z$, see, e.g, \cite{King}, \cite{Jan}. 
The subject  of the present paper
are closures of actions of infinite-dimensional groups. The first result of this kind was obtained by Nelson \cite{Nel},
who examined the closure of the  orthogonal group $\O(\infty)$ acting on the
infinite-dimensional space with Gaussian measure.  

Similar constructions related to non-singular transformations of measure spaces
 were used in ergodic theory,
see, e.g., \cite{Kre}, \cite{Dan}. 

\sm
 
The matter of author's interest are unitary representations
related to non-singular actions,  standard definitions are not relevant
to this topic. A nonsingular transformation $q$ of $Z$ produces a one-parametric
family of unitary operators in $L^2(Z)$, namely $f(z)\mapsto f(q(z))q'(z)^{1/2+is}$.
Therefore, a nonsingualar action of a group $G$ produces a one-parametric family of unitary representations of $G$.
So, a polymorphism must produce a one-parametric family of operators in $L^2(Z)$.
 Appropriate objects were introduced in   \cite{Ner-bist}.

\sm

{\bf\punct Polymorphisms.%
\label{ss:poly}
} 
We consider only Lebesgue  measure spaces, see \cite{Roh}, i.e., spaces, which 
are equivalent to a union of some interval in $\R$ and  at most countable collection of atoms having 
non-zero measure. 

For a space $Z=(Z,\zeta)$ with a continuous (non-atomic) measure  $\zeta$ we denote by
$\Ams(Z)$ or $\Ams(Z,\zeta)$ (an abbreviation of 'automorphisms of a measure space')  the group of all measure preserving transformations.
 There are 
two variants. If a measure $\zeta$ is finite, then without loss of generality
we can think that $Z$ is the segment $[0,1]$. If the measure is $\sigma$-finite,
then we can assume $Z=\R$. In the $\sigma$-finite case, we denote the group by
$\Aut_\infty(Z)$.
By $\Gms(Z)$ (or $\Gms(Z,\zeta)$) we denote the group of all non-singular transformations; in all cases, transformations are defined upto
sets of zero measure. Recall that a transformation $g:Z\to Z$ is {\it non-singular} if $g$ and $g^{-1}$ send sets of zero
measure to sets of zero measures. For a non-singular transformation $g$,
 we have a well-defined 
the {\it Radon--Nikodym} derivative $g'(z)$ 
determined from  the identity
$$
\zeta(g(A))=\int_A g'(x)\,d\zeta(z)\qquad \text{for all measurable $A\subset Z$.}
$$
We write the product in this group in the order $g_1\diamond g_2=g_2(g_1)$; equivalently, we
 consider the action of  $\Gms(Z)$  on $Z$ as a right action.

 Denote by $\R^\krest$ the multiplicative group of {\it positive} reals, denote by $t>0$ the coordinate on it.

\begin{definition} {\rm(\cite{Ner-bist})}
Let $(X,\xi)$, $(Y,\upsilon)$ be probabilistic Lebesgue spaces.
 A {\it polymorphism} 
$\frp: X\pol  Y$ is a measure on $X\times Y\times \R^\krest$
such that: 

\sm

{\rm 1.} the pushforward of $\frp$ under the projection 
$X\times Y\times \R^\krest\to X$ is the measure $\xi$;

\sm 

{\rm 2.} the pushforward of the measure $t\cdot \frp$ under the projection to $Y$ coincides with the measure
$\upsilon$.
\end{definition}

\sm 

We denote by $\Pol(X,Y)$ the space of all polymorphisms $X\pol Y$.

\sm 


We say that a polymorphism $\pi\in \Pol(X,Y)$ is {\it measure preserving} if
the measure $\pi$ is supported by the subset $X\times Y\times 1\subset  X\times Y\times \R^\krest$.
Objects of the previous subsection are measure preserving polymorphisms.

\sm 

{\sc Example.} Let $q\in \Gms(X,\xi)$. Consider the map 
$X\to X\times X\times \R^\krest$ defined by
$$
x\mapsto \bigl(x,q(x), q'(x)\bigr).
$$
Then the pushforward of the measure $\xi$ under this map is
a polymorphism $X\pol X$. If $q\in\Ams(X,\xi)$, then $q'=1$, and we get
a measure preserving polymorphism.
\hfill $\boxtimes$

\sm

{\sc The topology.}  Let $(M,\mu)$ be a  measure space,    $S\subset M$ be a measurable subset, let $\psi$ be a map from $M$ to a space $N$
sending $S$ to $T$.
We denote by $\mu\downarrow^S_T$ the pushforward of the restriction
of $\mu$ to $S$ under the map $\psi$ (in usages of this notation below,  $\psi$
is a projection map, so we omit $\psi$ and $N$ from the notation).

Conditions 1-2 in the definition of polymorphisms in this notation are: 
$$
\frp\dow^{X\times Y\times\R^\krest}_X=\xi, \qquad  \bigl(t\cdot \frp\bigr)\dow^{X\times Y\times\R^\krest}_Y=\upsilon.
$$

We say that a sequence  $\frp_j\in \Pol(X,Y)$ converges to $\frp$, 
if for any measurable subsets $A\subset M$, $B\subset N$ we have the following weak
convergences of measures%
\footnote{On weak convergences of measures on $\R$, see, e.g., \cite{Shi}, Sect.III.1-2.
Recall that the group $\R^\krest$ is isomorphic to the additive group of $\R$.} on $\R^\krest$:
\begin{align}
\frp_j\dow^{A\times B\times \R^\krest}_{\R^\krest}& \to \frp\dow^{A\times B\times \R^\krest}_{\R^\krest};
\label{eq:conv1}
\\
\bigl(t\cdot \frp_j\bigr)\dow^{A\times B\times \R^\krest}_{\R^\krest} &\to \bigl(t\cdot \frp\bigr)\dow^{A\times B\times \R^\krest}_{\R^\krest}.
\label{eq:conv2}
\end{align}

It is easy to show (see \cite{Ner-boundary}, Theorem 5.3) that {\it the group $\Gms(X)$ is dense in $\Pol(X,X)$}.

\sm 

{\sc Remark.} The spaces $\Pol(X,Y)$ are not compact (since a measure can disappear under a pass to a weak limit). I do not know, is it reasonable to compactify $\Pol(X,Y)$.
\hfill $\boxtimes$

\sm

{\sc Continuous polymorphisms.} See \cite{Ner-bist}.
Denote by $\frM^\triangledown(\R^\krest)$ the set of all positive finite Borel measures $\sigma$ on $\R^\krest$
such that $t\cdot \sigma$ is finite.
This set is a semigroup with respect to the convolution $*$ of measures on the group $\R^\krest$.
Indeed, 
$$
\int_{\R^\times} t\,d\sigma_1(t)*\sigma_2(t)=
\int_{\R^\times}\int_{\R^\times}pq\, d\sigma_1(p)\,\sigma_2(q)=
\int_{\R^\times}p\, d\sigma_1(p)\cdot \int_{\R^\times}q\, d\sigma_2(q).
$$

 A measurable function
$(x,y)\to s_{x,y}$
from $X\times Y$ to $\frM^\triangledown$ determines a  measure $\frs$
on $X\times Y\times \R^\krest$ in the following way. For measurable sets
$A\subset X$, $B\subset Y$, $C\subset \R^\krest$,
we assume
$$
\frs(A\times B \times C)=\int_{A\times B} s_{x,y}(C)\,d\xi(x)\,d\upsilon(y).
$$
We also must require that
$$
\int_{Y} s_{x,y}(\R^\krest)\,d\upsilon(y) =1,
\qquad
\int_X \int_{\R^\krest} t\, ds_{x,y}(t)\,d\xi(x) =1,
$$
under this condition we get an element of $\Pol(X,Y)$. We call such
polymorphisms {\it continuous}.

Let continuous polymorphisms
$\fru\in \Pol(X,Y)$, $\frv\in \Pol(Y,Z)$ be determined by functions
$(x,y)\mapsto u_{x,y}$, $(y,z)\mapsto v_{y,z}$.
Then their product 
$$\frw=\frv\diamond \fru \in \Pol(X,Z)$$
 is determined by the function
$$
w_{x,z}=\int_Y u_{x,y}*v_{y,z}\, d\upsilon(y).
$$ 
This is the formula for a product of integral operators, where the usual multiplication
is replaced by a convolution%
\footnote{If we consider only finite or countable measure spaces, then all polymorphisms are continuous.
Their product is reduced to a product of matrices over the ring of measures on $\R^\krest$
(or, more precisely, over the semiring $\frM^\triangledown$),
see \cite{Ner-boundary}, Subsect. 5.4.}.

\sm 

{\sc Product of polymorphisms.} It is easy to see, that continuous polymorphisms are dense
in $\Pol(X,Y)$. It can be shown that
{\it the multiplication $\diamond$ extends to a separately continuous map}
	$$
	\Pol(X,Y)\times \Pol(Y,Z)\to \Pol(X,Z). 
	$$
This operation is associative (see \cite{Ner-boundary}, Theorem 5.5, Theorem 5.9), 
i.e., for any Lebesgue spaces $X$, $Y$, $Z$, $U$ and for any 
$$
\frp\in\Pol(X,Y),\quad\frq\in \Pol(Y,Z),\quad \frr\in \Pol(Z,U),
$$
we have $(\frr\diamond \frq)\diamond \frp=\frr\diamond (\frq\diamond \frp)$.

 So, we get a category
 whose objects are Lebesgue probabilistic measure spaces and morphisms are polymorphisms. 

\sm

{\sc Remarks.}
a)
For other ways to define the multiplication of polymorphisms, see \cite{Ner-boundary}, however
  the definition by continuity and the definition by duality \eqref{eq:duality} given below are
most flexible.

\sm

b)  Informally,  polymorphisms are 'maps' $X\to Y$, which spread points of $X$ along $Y$; the Radon--Nikodym derivative also  spread.
 \hfill $\boxtimes$

\sm
{\sc The involution.}
Finally, we define  an involution on the category of polymorphisms.
For $\frp\in \Pol(X,Y)$, we consider its pushforward $\frp'$ under the map $(x,y,t)\mapsto (y,x,t^{-1})$  and  define $\frp^\star\in \Pol(Y,X)$ as the measure $t\cdot \frp'$. Then (see \cite{Ner-boundary}, Lemma 5.8)
for $\frp\in \Pol(X,Y)$, $\frq\in \Pol(Y,Z)$, we have
$$
(\frq\diamond \frp)^\star=\frp^\star \diamond \frq^\star.
$$

\bigskip

{\bf\punct Purposes of the paper.}
Descriptions of closures of some infinite-dimensional ({\it large}%
\footnote{The informal term 'large group' was proposed by A.~M.~Vershik. See a collection
of examples in Subsect. \ref{ss:examples}, only groups of types
A and C are infinite-dimensional in the formal sense; groups of other types 
are totally discontinuous.})
 groups in polymorphisms were obtained
in \cite{Ner-match}, \cite{Ner-gauss}, \cite{Ner-Whi}, the first two papers lead to unusual explicit formulas.

Quite often, large groups generate categories of double cosets ('trains').
 The aim of the present paper is to formulate two
general statements on this topic. In Theorem \ref{th:1}, we show that a non-singular action of a large
group $G$ extends to a functor from a category of double cosets 
 to the category of polymorphisms. The theorem is  simple (its proof was prepared in \cite{Ner-boundary}), however it is a general claim uniting quite different
 types of large groups (see Subsect. \ref{ss:examples}), such statements are relatively rare.
 In Theorem \ref{th:2}, we present a semi-constructive way to describe closures of large groups in polymorphisms.

\section{$\cM$-families of subgroups and trains of large groups}

\COUNTERS 

{\bf \punct Multiplicativity.%
\label{ss:multi}} Let $G$ be a separable topological group%
\footnote{We admit non-metrizable topologies, see Example A$_1$
in Subsect.
\ref{ss:examples}.}. 
Let $\rho$ be a unitary representation of $G$ in a {\it separable} Hilbert space $H$.
For any closed subgroup  $K\subset G$ we denote by $H^K$ the subspace
of $K$-fixed vectors, by $P^K$ the orthogonal projection to $H^K$. 
For two subgroups $K$, $L\subset G$ we define operators
$$
\wh \rho{\vphantom{\bigr|}}_{K,L}(g):H^L\to H^K
$$
by
\begin{equation}
\wh \rho\vph_{K,L}(g):=P^K \rho(g)\Bigr|_{H^L}.
\label{eq:wh}
\end{equation}
It is easy to see that 
$$
 p\in K,\, q\in L \quad \Longrightarrow\quad  \wh \rho\vph_{K,L}(pgq)=\wh \rho\vph_{K,L}(g).$$
So the operator $\wh \rho\vph_{K,L}(\cdot)$ actually depends on a double coset $\frg:=K\cdot g\cdot L$.
Denote the space of all double cosets by $K\backslash G/L$.  Next, we define a weaker equivalence relation
on $G$.

\begin{definition}
		 Reduced double cosets in $G$ with respect to subgroups $K$, $L$ are classes of the following equivalence:
	$$
	g\sim g' \,\Leftrightarrow \,
 \wh \rho\vph_{K,L}(g)=\wh \rho\vph_{K,L}(g')\,\text{for all unitary representations $\rho$ of $G$.}
	$$
	\end{definition}

Denote by $[K\backslash G/L]$ the set of all reduced double cosets. By the definition, 
$[K\backslash G/L]$ is a quotient space of $K\backslash G/L$.

\begin{definition}
A subgroup $K\subset G$ is an  $\cM$-subgroup in $G$ {\rm(}or a  pair $(G,K)$ satisfies  multiplicativity{\rm)} if for each $\frg$, $\frh\in [K\backslash G/K]$
there exists an element $\frg\cir \frh\in [K\backslash G/K]$ such that for any unitary representation
$\rho$ of $G$ we have
$$
\wh\rho\vph_{K,K}(\frg)\, \wh\rho\vph_{K,K}(\frh)=\wh\rho\vph_{K,K}(\frg\cir \frh).
$$
\end{definition}

Since a product of operators is associative, the $\cir$-product in $[K\backslash G/K]$
also is 
associative.

\sm

{\bf \punct  Remarks to the definition of $\cM$-subgroups.}
There is a big zoo of groups with $\cM$ subgroups, see below Subsect. \ref{ss:examples}.
But first of all, we will explain, where the multiplicativity cannot appear.

\sm

A) {\sc Hecke--Iwahori algebras $\frM(K\backslash G/K)$.}
 Let $G$ be a locally compact group and $K$ be a compact subgroup. Then  measures
on  $K\backslash G/K$ can be regarded as measures on $G$ invariant with respect to left and right shifts
by elements of $K$. So the space $\frM(K\backslash G/K)$ of all complex-valued Borel measures of bounded variation
 on $K\backslash G/K$
has a structure of an associative algebra.
 The product in this algebra  is the convolution $*$ of measures on $G$.
 
  Such algebras have been widely used
 in representation theory  since Gelfand theorem on spherical representations (1950).
 On one hand, they are used as  abstract tools; on the other hand, 
 some algebras of such types
 have nice structures and live with their own life (as well-known Hecke algebras, affine Hecke algebras, Yokonuma algebras).
 Except for few collections of examples, explicit description of such algebras are  difficult topics.
 
 For  $\frp\in K\backslash G/K$, we assign the unit measure $\delta_\frp$ supported by the point $\frp$.
 A convolution $\delta_\frp*\delta_\frq$
 determines an operation 
 $$
 K\backslash G/K\times K\backslash G/K\to \frM(K\backslash G/K),
 $$
 but $\delta_\frp*\delta_\frq$  is not a delta-measure, and we do not have
 a natural operation $K\backslash G/K\times K\backslash G/K\to K\backslash G/K$ (except a trivial case when
 $K$ is a normal subgroup, then $\frM(K\backslash G/K)$ is the convolution algebra of measures on the quotient
 group $G/K$). 
 
 \sm
 
B) {\sc 
  Our definition does not exclude trivial situations%
  \footnote{Our purpose is to formulate minimal conditions necessary to a proof of Theorem \ref{th:1}.}.} For instance, it
could happen that $G$ has no nontrivial unitary representations%
\footnote{Some examples are:
\newline	
	---  the group of operators $g$ in an infinite-dimensional Hilbert space such
that $g-1$ is a compact operator, see \cite{Nes0} (for a more transparent proof, see \cite{Pick}, Propositions 5.1.b and 7.1).
\newline
--- the group of all homomorphisms of the segment $[0,1]$, see \cite{Meg};
\newline 
 Considering overgroups of these groups,
one can produce numerous examples of groups without unitary representations.}.
Then  any $K\subset G$
is an $\cM$-subgroup, a semigroup $[K\backslash G/K]$ exists and consists of one element.
 
 \sm
 
C) {\sc For Lie connected groups, 
vectors fixed by noncompact subgroups are rare.} Let $G=\SL(n,\R)$ and $K\subset G$ be a noncompact closed subgroup.
  Then, for any nontrivial irreducible unitary representation $\rho$ in a Hilbert space $H$,
   we have%
     \footnote{Proof. The center $\cZ\simeq \Z/n \Z$ of $G$ consists of scalar  matrices $\theta$ such that $\theta^n=1$.
   	Consider the set $\rho(G)$ of all operators $\rho(g)$, where $g\in G$.
   	Denote by $\ov{\rho(G)}$   closure of $\rho(G)$ in the space of all operators equipped
   	with the weak  topology. According to Howe, Moore \cite{HM}, $\ov{\rho(G)}$ is a one-point compactification
   	of sume quotient of $G$ by a subgroup $U\subset \cZ$. The additional point 
   	is 0. Since $K$ is not compact $\ov{\rho(K)}=\rho(K)\cup 0$. If $K$ has a fixed vector $v$, then it is fixed by
   	all elements of $\ov{\rho(K)}$.  Therefore, $K$ has no fixed vectors.} $H^K=0$. Again we come to a one-point semigroup
  $[K\backslash G/K]$. 
  
  It is incorrect, that for  any connected Lie group, stabilizers of vectors in any unitary representation 
  are compact. But exceptions are semi-artificial (see a detailed discussion of this problem in \cite{Wan}) and are not interesting for our purposes. 
  
  For discrete groups, the situation  is not well-understood. A well-known interesting case is Example B$_1(2)$ in Subsect. \ref{ss:examples}.
  (also, Examples B$_1(n)$, D$_1(1)$). In fact, these  groups are 'large'. 


\sm

{\bf \punct $\boldsymbol\cM$-families and trains.} As above, let $G$ be a separable topological group
and $\rho$ be its unitary representation in a Hilbert space $H$. 
Let we have a subgroup $K=K_0\subset G$ and a family of subgroups  $K_\alpha\subset K$, where $\alpha$ ranges in some set $\cA$.
Assume that
for any pair $K_{\alpha}$, $K_{\alpha'}$ there is $K_{\alpha''}\subset K_{\alpha}\cap K_{\alpha'}$.
Let us denote
$$
H_\alpha:=H^{K_\alpha},\qquad P_\alpha:=P^{K_\alpha},
 \qquad \wh \rho_{\alpha,\beta}(\cdot):= \wh \rho\vph_{K_\alpha,K_\beta}(\cdot).
$$
Clearly,
$$
K_\alpha\supset K_{\alpha'}\,\, \Longrightarrow\,\, H_\alpha\subset H_{\alpha'},\,\, P_\alpha P_{\alpha'}= P_{\alpha}.
$$

\begin{definition}
A family $\{K_\alpha\}$ is an  $\cM$-family {\rm(}or a  family satisfying   multiplicativity),
 if for each $\alpha$, $\beta$, $\gamma\in \cA$ for any
$\frg\in [K_\alpha\backslash G/K_\beta]$,  $\frh\in [K_\beta\backslash G/K_\gamma]$
there is an element $\frg\cir \frh$ such that
for any unitary representation $\rho$ of $G$ we have
$$
\wh \rho_{\alpha,\beta}(\frg)\,\wh\rho_{\beta,\gamma}(\frh)=\wh\rho_{\alpha,\gamma}(\frg\cir\frh).
$$
\end{definition}

Since a product of operators is associative, $\cir$-product also is associative.

\begin{definition}
 The {\rm (}{\it train} $\Tr(G,\{K_\alpha\})$  is the category, whose
 objects are $\alpha\in \cA$, sets of morphisms are
$$
\Mor(\beta,\alpha)=[K_\alpha\backslash G/K_\beta],
$$
and a product is $\cir$-product.
\end{definition}

{\bf\punct Examples of $\boldsymbol\cM$-families.%
\label{ss:examples}}
We present (without proofs and with minimal explanations)
a broad variety of  examples of $\cM$-families $\{K_\alpha\}\subset G$. 

\sm

A. {\sc Infinite-dimensional real classical groups.} See, e.g., \cite{Olsh-motion},
\cite{Olsh-GB}, \cite{Pick}, \cite{Olsh-semi}, \cite{Ner-book}, Sect. IX.3-4, \cite{Ner-colligations}, 
\cite{Ner-spher}.

\sm 

--- $\mathrm{A}_1$.
 Let 
 $$G=\GL(\infty,\R)=\lim_{n\to\infty} \GL(n,\R)$$
  be the group of invertible real infinite matrices $g$
such that $g-1$ has only a finite number of nonzero elements.
We equip this group with the topology of inductive limit. Let $K=\O(\infty)$ be the subgroup consisting of orthogonal matrices.
Let $K_\alpha\subset K$ be the stabilizer  the first $\alpha$ elements of the basis. Then  $\{K\}_\alpha\subset G$
 is an $\cM$-family. Let us explain the product in $K_\alpha\backslash G/K_\alpha$. Consider the sequence
 \begin{equation}
 	\theta^{(\alpha)}_N=
 	\begin{pmatrix}
 		1_\alpha&0&0&0\\
 		0&0&1_N&0\\
 		0&1_N&0&0\\
 		0&0&0&1_\infty
 	\end{pmatrix},
 \end{equation}
 where $1_j$ denotes the unit matrix of size $j$. Then for any $g$, $h\in G$,
 the sequence of double cosets 
 $$\frq_{\vphantom{R^|}N}:=K_\alpha\cdot g \theta^{(\alpha)}_N h\cdot K_\alpha\,\in\, K_\alpha\backslash G/K_\alpha$$
  is eventually constant,
 its limit value is the product $\frg\cir \frh$. This is a general phenomenon%
 \footnote{Similar operations were known in spectral theory of non-self-adjoint operators 
 (see, e.g., \cite{Bro}, and system theory, see, e.g., \cite{Dym}, Chapter 19).} for 'representation theory of infinite-dimensional classical
 groups' in Olshanski sense \cite{Olsh-GB}.
 
 \sm  
 
--- $\mathrm{A}_2$.
 Let $G$ be the group of invertible operators
 in real $\ell^2$, which can be represented in the form $U(1+T)$, where $U$ is an orthogonal operator and $T$ is a Hilbert--Schmidt operator%
 \footnote{This group is the natural group of symmetries of a Gaussian measure
 on an infinite-dimensional space.}. Let $K$ be the group $\ov O(\infty)$ of all orthogonal operators.
 Subgroups $K_\alpha$ are defined in the same way. On the relations of $\mathrm{A}_2$ and $\mathrm{A}_1$,
 see the next subsection.
 
 \sm 
 
 --- $\mathrm{A}_3$. Let $G$, $K$  be as in the previous example. Let $\omega$
 be a subset of the set of odd numbers, $K_\omega$ be the subgroup in $\O(\infty)$
 fixing basis vector in $\ell^2$ with numbers $j\in \omega$.
 
 \sm
 
 {\sc B. Infinite symmetric groups.} See, e.g., \cite{Olsh-kiado}, \cite{Olsh-symm}, \cite{Olsh-semi}, \cite{Ner-umn}.
 
--- $\mathrm{B}_1(n)$.  Denote by $S_\infty$ the group of finitely supported permutations
 of $\N$,  by $\ov S_\infty$ the group of all permutations. Fix $n=1$, $2$, \dots.
 Let $G$ be the product of $n$ copies of $S_\infty$, where $n\ge 1$. 
 Let $K$ be the diagonal subgroup, so $K\simeq S_\infty$. Let $\{K_\alpha\}$ be point-wise stabilizers of
 sets $\{1, \dots,\alpha\}$. 
 
 \sm 
 
 --- $\mathrm{B}_2(n)$.
Consider the product of $n$ copies of $\ov S_\infty$. Let $G$ consist of tuples $(g_1,\dots,g_n)$ such that 
 $g_i g_j^{-1}\in S_\infty$ for all $i$, $j$. 
 Let $K\simeq\ov S_\infty$ be the diagonal subgroup,
  subgroups $\{K_\alpha\}$ are point-wise stabilizers of
 sets $\{1, \dots,\alpha\}$ as above. For the case $n=1$, i.e., the group $\ov S_\infty$, the classification
 of unitary representations was obtained by Lieberman \cite{Lie}, see Olshanski's comments in \cite{Olsh-kiado}.
  The case $n=2$ corresponds to a rich 'representation theory
 of infinite symmetric groups'
 initiated by E.~Thoma (A.~M.~Vershik, S.~V.~Kerov, A.~Wassermann, G.~I.~Olshanski, A.~Yu.~Okounkov, A.~M.~Borodin, et al.), see explanations in \cite{Olsh-symm}, \cite{Ner-umn}). 
 
  
  \sm
  
  {\sc C. Groups of transformations of measure spaces.} See \cite{Ner-bist}, \cite{Ner-book}, Sect. VIII.4, Chapter X, \cite{Ner-match}.

--- $\mathrm{C}_1$.   Let $Z$ be a probabilistic non-atomic measure space (so, we can set $Z=[0,1]$).
Let $G=\Gms(Z)$, $K=\Ams(Z)$. For any finite measurable partition $\frh$ of $Z$, we consider the subgroup $K_\frh\subset K$ consisting of transformations
preserving this partition.
In this case, we get the category of all polymorphisms of finite measure spaces
(and our definition of polymorphisms comes from this reasoning \cite{Ner-bist}).

 \sm
 
 --- $\mathrm{C}_2$.   Let $Z$ be a  non-atomic  measure space
 with  a $\sigma$-finite measure 
 (so, we can set $Z=\R$). Clearly,
 $\Gms([0,1])\simeq \Gms(\R)$. But we choose another $K=\Ams_\infty(Z)\simeq \Ams_\infty(\R)$. Let $\frh$ be a collection of disjoint measurable subsets
 $H_1$, \dots, $H_k$ of finite measure in  $Z$. We denote by $K_\frh$ the subgroup of $K$ preserving all subsets 
 $H_j$. Polymorphisms related to this $\cM$-family are discussed
 in Subsect. \ref{ss:final}
  
  \sm 
  
--- $\mathrm{C}_3$.   Let $Z$ be a  non-atomic measure space
with a $\sigma$-finite measure.  Let $G$ be the group
of nonsingular transformations $g$ of $Z$ such that%
\footnote{This group is the natural group of transformations of a  Poisson
point processes, see \cite{Ner-book}, Sect. X.4.} $(g'(z)-1)\in L^1(Z)$. The subgroup $K$
and subgroups $K_\frh$ are defined as in  $\mathrm{C}_2$. See \cite{Ner-book}, Sect. X.4, \cite{Ner-match}.

\sm  
  
  {\sc D. Infinite-dimensional matrix groups over finite fields.}
  See \cite{Olsh-semi}, \cite{Ner-finite}, \cite{Tsa}. Let $\mathbb{F}_q$ be a finite field.
  
  ---  $\mathrm{D}_1$. Let $G$ be the group of infinite matrices $g$ over  $\mathbb{F}_q$ such that $g-1$
  has only a finite number of nonzero matrix elements. Let $K_\alpha$ consist of matrices of the form
  $\begin{pmatrix}
  	1_\alpha&0\\ 0& *
  \end{pmatrix}$.

\sm  
  
---   $\mathrm{D}_2$.
Let $V$ be the space of sequences $(x_1, x_2,\dots)$, where $x_j\in \mathbb{F}$, such that all but a finite number of elements are
zero. Let $G$ be the group of all linear transformations of $V$. Let $K=G$ and $K_\alpha$ be the stabilizer of the first $\alpha$ elements of the basis.

\sm 

---   $\mathrm{D}_3$. Consider the group $G$ of invertible matrices over $\mathbb{F}_q$, which have  finite number of nonzero elements in each
column and each row. Denote by $e_j$ the standard basis in the spaces of vector-columns, by $f_i$ the standard basis
in the space of vectors-row. Let $A$, $B$ be finite subsets in $\N$. Denote by $K_{A,B}\subset G$
the subgroup fixing all vectors $e_a$, where $a\in A$, and $f_b$, where $b\in B$.

\sm

--- $\mathrm{D}_4$. Let $W$ be the space of two-side sequences $\{x_j\}$,
where $j\in \Z$, $x_j\in \mathbb{F}_q$, such that $x_j=0$ for sufficiently
large $j$. The $W$ is a unique separable locally compact linear space 
over $\mathbb{F}_q$ that is not discrete and not compact. Denote by $W_\alpha$ the subspace
consisting of sequences such that $x_j=0$ for $j>\alpha$. Denote by
$G$ the group of all linear transformations  of $W$. For two
subspaces $W_\alpha\supset W_\beta$ denote by $K_{\alpha,\beta}\subset G$
the subgroup preserving the flag $W_\alpha\supset W_\beta$ and inducing the trivial linear transformations on the quotient $W_\alpha/W_\beta$.

\sm

{\sc E. Infinite-dimensional matrix groups over $p$-adic fields.}
See., e.g., \cite{Ner-adic}.

--- $\mathrm{E}_1$. Let $G$ be the group of invertible infinite matrices over 
$p$-adic integers, which have only finite number of elements in each row.
Let $K_\alpha$ be the subgroup fixing first $\alpha$ basis elements.

\sm 

 {\sc F. Groups of automorphisms and spheromorphisms
 of non-locally finite trees and $\R$-trees.}
 See, e.g.,
\cite{Olsh-new}, \cite{Ner-book}, Sect. VIII.6 (Remarks), \cite{Ner-urysohn}.
 
--- $\mathrm{F}_1$. 
  Consider a simplicial tree $T$  such that each vertex is contained
in a countable family of edges. Let $G$ be the group of its automorphisms, $K$ be a stabilizer of a vertex, say $w$.
Subgroups $K_\alpha$ are stabilizers of finite subtrees containing $w$. 

\sm

 {\sc G. Oligomorphic groups%
 \footnote{A closed subgroup in $\ov S_\infty$ is {\it oligomorphic} if for each $m$
 it has a finite number of orbits on $\N^m$.}. } See \cite{Ner-universal}, in fact this is consequence of Tsankov \cite{Tsa}.

--- $\mathrm{G}_0$. See  examples  $\mathrm{D}_2$, $\mathrm{D}_3$, $\mathrm{F}_1$.

\sm 

--- $\mathrm{G}_1$. 
Let $G$ be the group of automorphisms of the Rado graph (see, e.g., \cite{Cam}). Subgroups $K_I$ are stabilizers of
finite subgraphs $I$. 

\sm

--- $\mathrm{G}_2$. Let $G$ be the group of all order preserving bijections of $\Q$. Subgroups $K_A$
are stabilizers of finite subsets in $\Q$.
  
 \sm
 
{\bf\punct Remarks on examples of $\boldsymbol\cM$-families.} 

\sm
 
 {\sc Examples of Ismagilov.}   The first examples of $\cM$-subgroups were discovered by Ismagilov \cite{Ism69} (1967),
  \cite{Ism73} (1973). Namely, he considered the following situations:
  
  \sm 
  
--- $\mathrm{F}_2$.  Let $\K_1$ be a complete  non-Archimedean discrete valuation field 
$\K_1$ with a countable field of residues, let $\mathbb{O}_1$
be the ring of integers of $\K_1$. We set $G_1=\SL(2,\K_1)$, $K_1=\SL(2,\mathbb{O}_1)$.

\sm 

--- $\mathrm{F}_3$.
  Let $\K_2$ be a complete  non-Archimedean valuation field 
  with a non-discrete  group of valuation. Denote by 
 $\mathbb{O}_2$ the ring of integers.   We set $G_2=\SL(2,\K_2)$, 
 $K_2=\SL(2,\mathbb{O}_2)$.

\sm 

The group $G_1=\SL(2,\K_1)$ acts on a non-locally finite simplicial tree, and $G_2=\SL(2,\K_2)$
on an $\R$-tree, the subgroups  $K_j=\SL(2,\mathbb{O}_j)$ are stabilizers of vertices, see, e.g. \cite{Chi}.
Strangely enough, that
Ismagilov operated with trees, $\R$-trees, and spherical functions on trees  but he
  didn't see trees. %
The paper \cite{Ism69} was an important standpoint for Olshanski, see \cite{Olsh-motion},
 \cite{Olsh-new}. 
 
 \sm 

  {\sc On the reduction of double coset spaces.}  In many cases, unitary representations separate double cosets,
  so $[K\backslash G/L]=K\backslash G/L$. A situation, when $K\backslash G/L$ is not a Hausdorff
  topological space, also is usual (as  cases $\mathrm{C}_1$--$\mathrm{C}_3$, see  \cite{Ner-canonical}).
   This is one of reasons for a reduction. However,
  Olshanski \cite{Olsh-semi} observed that a reduction can appear even for finite spaces $K\backslash G/K$
  with discrete topology, namely for the case $\mathrm{D}_1$.
  
\sm 

{\sc Admissible representations.} 
For Theorem \ref{th:2} below, we need more details and more definitions.

\begin{definition}
Let	$\rho$ be a unitary representation of $G$ in $H$. Denote 
$$H^\circ:=\bigcup_{\alpha\in\cA}
 H_\alpha.
 $$
 We say that $\rho$ is:
 
 \sm 
 
 ---   admissible if $H^\circ$ is dense in $H$;

\sm

---  essentially admissible if there are no proper subrepresentations
in $H$ containing $H^\circ$.
\end{definition} 

Let $\rho$ be an arbitrary unitary representation of $G$. We consider its minimal subrepresentation
$H_{\mathrm{adm}}$ containing $H^\circ$ and its orthocomplement $H_{\mathrm{adm}}^\bot$.
Then we get an essentially admissible representation of $G$ in $H_{\mathrm{adm}}$
and a representation without $K_\alpha$-fixed vectors in $H_{\mathrm{adm}}^\bot$. 
Clearly, the train construction can be interesting only for
the  essentially admissible summand.

Next, 
the group from example $\mathrm{A}_1$ has not type I, the set of its irreducible
unitary representations is incredible, and there is no chance to understand it.
The group from example $\mathrm{A}_2$ is type I, classification of its irreducible
representations was conjectured by Olshanski in \cite{Olsh-GB}
(see, also, an exposition in \cite{Ner-book}, Sect. IX.5), the conjecture is doubtless
(but  it has not been proved until now). In any case, admissible representations of 
groups $\mathrm{A}_1$ and $\mathrm{A}_2$ coincide (see, e.g., Pickrell \cite{Pick}, Proposition 5.1.c, this also follows from \cite{Nes00}).
The same holds for $\mathrm{B}_1$ and $\mathrm{B}_2$, $\mathrm{D}_1$ and $\mathrm{D}_3$.


In Examples $\mathrm{A}_2$, $\mathrm{B}_2(n)$, $\mathrm{C}_1$--$\mathrm{C}_3$, $\mathrm{D}_2$--$\mathrm{D}_4$,
$\mathrm{E}_1$, $\mathrm{F}_1$, $\mathrm{G}_1$--$\mathrm{G}_2$, all continuous representations
are admissible%
\footnote{In fact, admissibility is a property of the subgroup $K$ (and not of its overgroup $G$.}. 

\begin{definition}
	\label{def:admis}
We say that an $\cM$-family is exhaustive if: 

\sm 

{\rm 1)} any essentially admissible representation of $G$ is admissible;

\sm

{\rm 2)} for each $K_\alpha$ there is a sequence $\theta^{(\alpha)}_j\in K_\alpha$
such that for any unitary  representation $\rho$ the sequence 
$\rho\bigl(\theta^{(\alpha)}_j)$ converges
to $P_\alpha$ in the weak operator topology. 

\sm 

{\rm 3)} there is a sequence of subgroups $K_{\gamma_1}\supset K_{\gamma_2}\supset \dots$
such that 
for any unitary admissible representation of $G$ the subspace
 $\cup_j H_{\gamma_j}$ is dense%
 \footnote{It sufficient to assume that for any $K_\alpha$
 	 there is $N$ such that $K_{\gamma_N}\subset K_\alpha$. This holds automatically in all examples
  except $\mathrm{C}_1$--$\mathrm{C}_3$.}.
\end{definition}

In all examples known to the author,  conditions 2) and 3) hold automatically.
In  example $\mathrm{A}_3$ from our list, the condition 1) fails. But this $\cM$-family
of subgroups can be easily extended to an exhaustive  $\cM$-family (we say that $\omega$ is a subset of $\N$ with an infinite complement,
the new $\cM$-family includes the $\cM$-family from $\mathrm{A}_2$.

\sm

{\sc Change of topologies and definitions of $\cir$-products.}
For all examples $(G,\{K_\alpha\})$ of $\cM$-families known to the author,
we can pass to a completion $G^\bullet$ of $G$ such that 
$$
\Bigl\{\text{admissible representations of $G$}\Bigr\}\simeq \Bigl\{\text{continuous representations of
	$G^\bullet$}\Bigr\},
$$
as $\mathrm{A}_1\to \mathrm{A}_2$, $\mathrm{B}_1\to \mathrm{B}_2$, $\mathrm{D}_1\to \mathrm{D}_3$. Generally, this completion is not unique, in all examples
known to the author there is a (may be, non-unique) completion with a Polish topology%
\footnote{A topological space is {\it Polish} if it is homeomorphic to a complete separable metric space.}.

On the other hand, in example $\mathrm{A}_1$ we presented a simple constructive definition
of products of double cosets. Usually, we can pass to groups with more strong separable
(may be, non-Polish)
topologies with the same set of admissible representations, for which this approach is effective.

\section{Actions of trains by polymorphisms}

\COUNTERS

Here we discuss the following situation.
Let $G$ be a separable topological group, $\{K_\alpha\}$ be an $\cM$-family of subgroups.
Consider an action of $G$ by nonsingular transformations $\sigma(g)$ of a Lebesgue probabilistic space
$(Z, \zeta)$, i.e., a continuous homomorphism%
\footnote{We can define a topology on $\Gms(Z)$ as the topology induced from $\Pol(Z,Z)\subset \Gms(Z)$.
On the other hand, we can fix $1<p<\infty$ and consider the group $\mathrm{Isom}(L^p)$
 of all isometries of $L^p(Z)$ equipped
with a weak or strong operator topology. Next, we fix $s\in \R$ and embed $\Gms(Z)$ to the group
$\mathrm{Isom}(L^p)$ by the formula \eqref{eq:operator} with $r=1/p$ and restrict the topology
of $\mathrm{Isom}(L^p)$ to $\Gms(Z)$. Then all such induced topologies
coincide
with the topology induced from $\Pol(Z,Z)$. A formal proof is contained in \cite{Ner-canonical},
but this is superfluous. The statement 'any homomorphism between two Polish groups is continuous'
is compatible with Zermelo--Fraenkel system + the Axiom of dependent choice, see, e.g., \cite{Wri}. 
 In particular, if we define constructively two structures of a Polish topological group
on the same abstract group, then these topologies will inevitably coincide.}
 from $G\to \Gms(Z)$.
Assume that $K$ acts by measure preserving transformations.

We wish to show that such action automatically produces an action of 
the train $\Tr(G,\{K_\alpha\})$ by polymorphisms.

\sm

{\bf \punct Sigma-algebras of $K_\alpha$-invariant sets.}
Let a group $H$ act by measure preserving transformations
 on a Lebesgue space $(Y,\upsilon)$.
We say, that a measurable set $A\subset Y$ is $H$-{\it invariant} if
for any $h\in H$  the symmetric difference $A\bigtriangleup h(A)$
has zero measure. Clearly, $H$-invariant sets form a $\sigma$-algebra.

Consider an action of $G$ as above. We claim that

\sm 

--- {\it for each $\alpha$ there exists  a well-defined quotient, say $(Z_\alpha, \zeta_\alpha)$, of
$(Z,\zeta)$ by the action of the group $K_\alpha$}.

\sm 

 Namely, consider
the sigma-algebra $\Sigma_\alpha$ of all $K_\alpha$-invariant measurable subsets
 in $Z$.
Also consider the space $H_\alpha$ consisting of $K_\alpha$-fixed functions in
$L^2(Z,\zeta)$. For each $A\in \Sigma_\alpha$ we assign its indicator function
$I_A(z)\in L^2(Z)^{K_\alpha}$. Since $L^2$ is separable, we can choose a dense countable subset $I_{A_j}$
in the space of all functions $I_A$, where $A$ ranges in $\Sigma_\alpha$.
Consider the sigma-algebra $\Sigma'_\alpha$ generated by $A_j$. According
Rohlin \cite{Roh}, it determines a measurable partition of $Z$.
Denote by $Z_\alpha$ the quotient space with respect to  this partition. We have a 
canonical map $Z\to Z_\alpha$ and therefore we get a measure, say $\zeta_\alpha$,
on $Z_\alpha$. \hfill $\square$

Let $(X,\xi)$, $(Y,\upsilon)$  be probabilistic  measure spaces, let $\pi:X\to Y$
be a measurable map such that for each measurable $B\subset Y$ we have
$\xi(\pi^{-1}(B))= \upsilon(B)$.
We define a polymorphism $\frl[\pi]:X\pol Y$ as the image of $\xi$ under the map
 $X\to X\times Y\times \R^\krest$ given by
$x\mapsto (x,\pi(x),1)$, see \cite{Ner-boundary}, Subsect. 3.10. 

Since we have a canonical map $Z\to Z_\alpha$, we have a canonical polymorphism
$\frl_\alpha:Z\pol Z_\alpha$.

{\bf \punct The statement of the paper.}

\begin{theorem}
\label{th:1}
Let $G$ be a separable topological group, $\{K_\alpha\}_{\alpha\in\cA}$ be an $\cM$-family of subgroups.
Consider an action of $G$ by nonsingular transformations $\sigma(g)$ of a Lebesgue probabilistic space
$Z$. Let $K$ act by measure preserving transformations. 
For any $\alpha$, $\beta\in\cA$ we define
$$
\cS_{\alpha,\beta}(g)=\frl_\alpha^\star\diamond \sigma(g)\diamond \frl_\beta\in \Pol(Z_\beta,Z_\alpha).
$$
Then $\cS$ is a functor from the train $\Tr(G,\{K_\alpha\}_{\alpha\in \cA})$ to the category of polymorphisms.
Namely,
$\cS_{\alpha,\beta}(g)$ depends only on a reduced double coset 
$\in[K_\alpha\backslash G/K_\beta]$ containing $g$
and  for any 
$\frg\in [K_\alpha\backslash G/K_\beta]$, $\frh\in [K_\beta\backslash G/K_\gamma]$
we have
$$
\cS_{\alpha,\beta}(\frg)\diamond \cS_{\beta,\gamma}(\frh)=
\cS_{\alpha,\gamma}(\frg\cir\frh).
$$
\end{theorem}

We need some preliminaries.

\sm

{\bf \punct Mellin--Markov transform of polymorphisms.}
For a Lebesgue probabilistic measure space $Z$ denote by $L^{\infty-}(Z)$ the space 
of bounded measurable functions on $Z$ equipped with the following convergence:
$\phi_j\to \phi$ if essential supremums of $|\phi_j|$ are uniformly bounded and the sequence $\phi_j$ converges to
$\phi$ in measure. 

Let $r+is\in\C$ range in the strip
$\Pi: 0\le r\le 1$. 
For $\frp\in \Pol(X,Y)$ and $r+is\in \Pi$,
we define the bilinear form $B_{r+is}:L^{\infty-}(X)\times L^{\infty-}(Y)\to \C$ by
$$B_{r+is}[\frp](\phi,\psi)=
\iint\limits_{X\times Y \times \R^\krest}
\phi(x)\, \psi(y)\, t^{r+is} \,d\frp(x,y,t).
$$
We define the {\it Mellin--Markov transform of $\frp$}
as an operator-valued function $r+is\mapsto T_{r+is}(\frp)$,
$$
T_{r+is}(\frp):L^{1/r}(Y)\to L^{1/r}(X)
$$
determined by the following duality:
\begin{equation}
\int\limits_X \phi(x)\cdot \bigl(T_{r+is}(\frp) \psi(x)\bigr)\,d\xi(x)=B_{r+is}[\frp](\phi,\psi)
\end{equation}
for any $\phi\in L^{\infty-}(X)$, $\psi\in L^{\infty-}(Y)$.
Then $\|T_{r+is}(\frp)\|_{L^{1/r}}\le 1$ for $o\le r\le 1$.
Also, $T_{0+is}(\frp)$ is a continuous operator $L^{\infty-}(Y)\to L^{\infty-}(X)$,
see \cite{Ner-boundary}, Theorem 6.3.

\sm

{\sc Example.} Let $g\in \Gms(Z)$ be considered as a polymorphism. Then
\begin{equation*}
\label{eq:operator}
\qquad\qquad\qquad\qquad\qquad
T_{r+is}(g)\,\phi(z)=\phi(g(z))\,g'(z)^{r+is}. 
\qquad\qquad\qquad\qquad\qquad\boxtimes
\end{equation*}

\sm 

A polymorphism $\frp$ is uniquely determined by its Mellin--Markov transform (\cite{Ner-boundary}, Theorem 6.12).
For any $\frp\in\Pol(M,N)$, $\frq\in\Pol(N,K)$ we have (see \cite{Ner-boundary}, Theorem 6.14)
\begin{equation}
T_{r+is}(\frq)\, T_{r+is}(\frp)=T_{r+is}(\frq\diamond \frp).
\label{eq:duality}
\end{equation}

This transform is holomorphic in the variable $r+is$ in the following sense:
for any bounded measurable functions $\phi$, $\psi$, the matrix element
$r+is\mapsto B_{r+is}(\phi,\psi)$ is holomorphic in the strip $0<r<1$ and continuous in the closed strip $0\le r\le 1$.
So, $T_{r+is}(\frp)$ is an analytic family of operators in the strip $\Pi$ in the sense of Stein \cite{Ste}. 

Pointwise convergence of the  bilinear forms 
$$B_{r+is}[\frp_j](\phi,\psi)\to B_{r+is}[\frp](\phi,\psi)$$
in the closed strip $\Pi$
is equivalent to the convergence of polymorphisms $\frp_j\to\frp$ (\cite{Ner-boundary}, Theorem 6.14).

\sm

{\bf \punct Proof of Theorem \ref{th:1}.}
For $\alpha\in \cA$ we denote by $L^2(Z)_\alpha$ the subspace consisting of $K_\alpha$-fixed functions in $L^2(Z)$.
Denote
\begin{equation}
	\frm_\alpha:=\frl_\alpha^\star \diamond \frl_\alpha \in \Pol(Z,Z).
\end{equation}
Clearly, $\frm_\alpha^\star=\frm_\alpha$.
Notice that $ \frl_\alpha \diamond  \frl_\alpha^\star$ is the unit in $\Pol(Z_\alpha,Z_\alpha)$.
Therefore, $\frm_\alpha^2=\frm_\alpha$.

\sm

Firstly, let us describe the operators $T_{1/2+is}[\frl_\alpha]$, $T_{1/2+is}[\frl_\alpha^\star]$, $T_{1/2+is}[\frm_\alpha]$
(see, e.g, \cite{Ner-boundary}, Subsect. 3.9). The polymorphsims $\frl_\alpha$, $\frl_\alpha^\star$, $\frm_\alpha$
are measure preserving,  so their Mellin--Markov transforms do not depend on $s$.

\sm

--- The operator $T_{1/2+is}[ \frm_\alpha]$
is the operator of orthogonal projection $P_\alpha:L^2(Z)\to L^2(Z)_\alpha$ or, equivalently, the operator of conditional
expectation, see 
\cite{Ner-boundary}, Subsect. 3.10.

\sm 

--- The operator $T_{1/2+is}[\frl_\alpha^\star]:L^2(Z_\alpha)\to L^2(Z)$
is an  isometric embedding, its image coincides with $L^2(Z)_\alpha$.

\sm 

--- The operator
$$
T_{1/2+is}[\frl_\alpha]\Bigr|_{L^2(Z)_\alpha}:\,L^2(Z)_\alpha\to L^2(Z_\alpha)
$$
is unitary, and
 $$
 T_{1/2+is}[\frl_\alpha]\Bigr|_{(L^2(Z)_\alpha)^\bot}=0.
 $$

\sm 

 Secondly, we slightly reformulate the $\cM$-property.
In notation of Subsect. \ref{ss:multi}, we define operators $\wt\rho_{\alpha,\beta}(g)$ in $H$ by
$$
\wt \rho_{\alpha,\beta}(g)=P_\alpha\, \rho(g)\, P_\beta.
$$
These operators have the following block form:
$$
\wt\rho_{\alpha,\beta}(g)=\begin{pmatrix} \wh \rho_{\alpha,\beta}(g)&0\\0&0\end{pmatrix}:
H_\beta\oplus H_\beta^\bot\to H_\alpha\oplus H_\alpha^\bot,
$$
where operators $\wh\rho(\cdot)$ are defined by \eqref{eq:wh}.
Clearly, these operators depend only on double cosets, and
the $\cM$-property can be written in the form
\begin{equation}
\wt\rho_{\alpha,\beta}(\frg)\,\wt\rho_{\beta,\gamma}(\frh)=\wt\rho_{\alpha,\gamma}(\frg\cir\frh).
\label{eq:1}
\end{equation}

In the same way, we define polymorphisms
$$
\wt \cS_{\alpha,\beta}(g):= \frm_\alpha \diamond \sigma(g)\diamond \frm_\beta \in \Pol(Z,Z).
$$
The conclusion of the theorem can be reformulated in the form
\begin{equation}
	\text{[\,{\huge ?}\,]} 
	\qquad
\wt \cS_{\alpha,\beta}(\frg)\diamond \wt \cS_{\beta,\gamma}(\frh)=\wt \cS_{\alpha,\gamma}(\frg\cir \frh)
 .
\label{eq:2}
\end{equation}

Applying the Mellin-Markov transform to both sides of the hypothetical equality \eqref{eq:2} we get
an equivalent hypothetical equality  
$$
	\text{[\,{\huge ?}\,]} 
\qquad
T_{r+is}\bigl(\wt \cS_{\alpha,\beta}(\frg)\diamond \wt \cS_{\beta,\gamma}(\frh)\bigr)=
T_{r+is}\bigl(\wt \cS_{\alpha,\gamma}(\frg\cir \frh)\bigr).
$$
Expressions in both sides are holomorphic in the sense of Stein. We rewrite this in the form
\begin{multline}
		\text{[\,{\huge ?}\,]} 
	\qquad
T_{r+is}(\frm_\alpha)\, T_{r+is}(g) \, T_{r+is}(\frm_\beta)\,
T_{r+is}(h) T_{r+is}(\frm_\gamma)
=\\=
 T_{r+is}(\frm_\alpha)\,  T_{r+is}(g\cir h)\, T_{r+is}(\frm_\gamma),
\label{eq:3}
\end{multline}
where operators $T_{r+is}(g)$ are given by \eqref{eq:operator}, and $g\cir h$ denotes an arbitrary representative of the product of double cosets.
We set  $r=1/2$. Keeping in  mind that $T_{r+is}(\frm_\alpha)=P_\alpha$
 we come to the following  hypothetical	 equality of operators in $L^2(Z)$:
$$
	\text{[\,{\huge ?}\,]} 
\qquad
P_\alpha\, T_{1/2+is}(g)\, P_\beta\, T_{1/2+is}(h) P_\gamma= P_\alpha\, T_{1/2+is}(g\cir h)\, P_\gamma.
$$
But this is the identity \eqref{eq:1} for the representations $T_{1/2+is}\cir \sigma$ of $G$. So \eqref{eq:3} holds on the line $r=1/2$. By holomorphy of
both sides of \eqref{eq:3}, this holds in the whole strip $\Pi$. This implies  \eqref{eq:2}.

\sm

{\bf \punct Closures of groups in polymorphisms.} 
Let $G$ be a separable topological group, $\{K_\alpha\}_{\alpha\in\cA}$ be an exhaustive $\cM$-family of subgroups.
Consider an action of $G$ by nonsingular transformations $\sigma(g)$ of a Lebesgue probabilistic space
$Z$. Let $K$ act by measure preserving transformations.   
Let $K_{\gamma_k}$ is a sequence of subgroups as in Definition \ref{def:admis}.

\begin{theorem}
	\label{th:2} 
	Assume that the representation $T_{1/2+is}$
	of $G$ 
	in $L^2(Z)$ is admissible%
	\footnote{The admissibility is defined in terms of $K$, therefore this condition does not depend on
		$s$.}	
A polymorphism $\frr\in \Pol(Z,Z)$ is contained in the closure of $G$ if and only if
for any $\gamma_k$ the element
$\frm_{\gamma_k}\diamond \frr \diamond\frm_{\gamma_k}$ 
is contained in the closure of $G$ for all $k$.
\end{theorem}

{\sc Proof.} The statement $\Rightarrow$.
For a fixed $\gamma_k$ consider a sequence $\theta^{(\gamma_k)}_j$ as in Definition \ref{def:admis}. 
Then $T_{1/2+is}\bigl(\theta^{(\gamma_k)}_j \bigr)$ converges to 
$P_{\gamma_j}$. Equivalently, $\theta^{(\gamma_k)}_j\in\Ams(Z)$ converges to $\frm_{\gamma_j}$ in
(measure preserving) polymorphisms. So,  $\frm_{\gamma_j} \diamond \frr \diamond\frm_{\gamma_j}$ 
is contained the closure.

\sm 

 The statement $\Leftarrow$.
 Let us show that $\frr$ is the limit 
of $\frm_{\gamma_j}\diamond \frr\diamond \frm_{\gamma_j}$. It is sufficient to prove
 the pointwise convergence of corresponding bilinear forms $B_{r+is}$,
i.e.,
$$
B_{r+is}[\frm_{\gamma_k} \diamond\frr\diamond \frm_{\gamma_k}](\phi,\psi)\to B_{r+is}[\frr](\phi,\psi).
$$
We transform this expression as 
\begin{multline*}
B_{r+is}[\frm_{\gamma_k} \diamond \frr \diamond \frm_{\gamma_k}](\phi,\psi)
=\int_Z \phi(z)\cdot  T_{r+is} (\frm_{\gamma_k} \diamond \frr \diamond \frm_{\gamma_k})\psi(z)\,d\zeta(z)
=\\=
\int_Z \phi(z)\cdot   (P_{\gamma_k}  T_{r+is}(\frr) P_{\gamma_k})\psi(z)\,d\zeta(z)
=\\= \int_Z (P_{\gamma_k}\phi)(z)\cdot     T_{r+is}(\frr) (P_{\gamma_k}\psi(z))\,d\zeta(z)
=
B_{r+is}[\frr](P_{\gamma_j} \phi,P_{\gamma_j}\psi).
\end{multline*}
So we must verify the convergence
\begin{equation}
B_{r+is}[\frr](P_{\gamma_k} \phi,P_{\gamma_k}\psi) \to B_{r+is}[\frr](\phi,\psi),
\label{eq:hide}
\end{equation}
where $\phi$, $\psi$ are bounded functions.

Recall that $P_{\gamma_k}$ are operators of conditional expectations.
Therefore, the sequences $P_{\gamma_k} \phi$, $P_{\gamma_k} \psi$
are bounded martingales, so they converge a.s and in the $L^1$-sense (see, e.g., \cite{Shi}, Theorems VIII.4.1-2). 
By the density of $\cup L^2(Z)_{\gamma_k}$, our sequences converge to 
$\phi$, $\psi$ in $L^2$, and therefore we have a.s. convergences $P_{\gamma_k} \phi\to \phi$, $P_{\gamma_k} \psi\to\psi$.

Next, 
\begin{equation*}
B_{r+is}[\frr](P_{\gamma_k} \phi,P_{\gamma_k}\psi)=\int_{Z\times Z\times \R^\krest}
\Bigl(P_{\gamma_k} \phi(z_1)\, P_{\gamma_k} \psi(z_2)\, t^{is} \Bigr)\cdot t^r \,d\frr(z_1,z_2,t).
\end{equation*}
By the definition of polymorphisms, a measure $t^r \frr(z_1,z_2,t)$ is finite
(since $t^r\le \max \{1,t\}$). For fixed $\phi$, $\psi$, the expressions
in the big brackets are uniformly bounded. Applying  Lebesgue's dominated convergence theorem
we get \eqref{eq:hide}. \hfill $\square$



\sm

{\bf \punct Some questions.%
\label{ss:final}} 
1)  The main question is examinations of explicit non-singular (and measure-preserving) actions
and search of  explicit formulas for actions of trains and closures of actions.
We do not try to review existing constructions of actions of large groups on measure spaces (and I do not know such surveys), some important examples are inverse limits of symmetric groups \cite{Ker} and inverse limits of classical groups
and compact symmetric spaces \cite{Pick0}, \cite{Ner-hua}. This paper is abstractionist,
 and we add some abstractionist remarks. 
 
 \sm 
 
2) Apparently, there is a version of Theorem \ref{th:1} for spaces $(X,\xi)$ with $\sigma$-finite measures.
Namely, let we have a space $(Z,\xi)$ with a non-atomic $\sigma$-finite measure, we have a group
$G$ with an $\cM$-family of subgroups $\{K_\alpha\}$. Let we have a non-singular action of $G$
on $Z$, let  $K$ act by measure preserving transformations.

 In this case, a definition of polymorphisms relevant to
 unitary representations must keep in mind that a measure can go to infinity, and a measure can come from infinity.
 Moreover, a  measure can depart and arrive  with different Radon--Nikodym derivatives. The technics of
 \cite{Ner-bist} leads to the following definition.
 
 

\begin{definition} 
Let $(X,\xi)$, $(Y,\upsilon)$ be spaces with  $\sigma$-finite Lebesgue measures. We say that an element
of $\Pol_\infty(X,Y)$ is a triple $(\frr, \frr^-, \frr^+)$, where $\frr$ is a measure on $X\times Y\times \R^\krest$,
$\frr_-$ is a measure on $X\times \R^{\krest}$, $\frr_+$ is a measure on $Y\times \R^\krest$. This 
triple must satisfy two conditions:
\begin{align*}
 \frr\dow^{X\times Y\times \R^\krest}_X+ 
 \frr^-\dow^{X\times \R^\krest}_X=\xi,\qquad 
\bigl(t\cdot  \frr\bigr)\dow^{X\times Y\times \R^\krest}_Y+ 
\bigl(t\cdot \frr^+\bigr)^{Y\times \R^\krest}_Y=\upsilon.
\end{align*} 
\end{definition}

We say, that a sequence $(\frr_j, \frr^-_j, \frr^+_j)$ converges to $(\frr, \frr^-, \frr^+)$
if for any measurable sets $A\subset X$, $B\subset Y$ of finite measure we have the following
weak convergences of measures on $\R^\krest$ (the first two lines correspond to conditions
\eqref{eq:conv1}--\eqref{eq:conv2}):
\begin{align*}
\frr_j\dow^{A\times B\times \R^\krest}_{\R^\times }&\to
\frr\dow^{A\times B\times \R^\krest}_{\R^\times };
\\
\bigl(t\cdot\frr_j\bigr)\dow^{A\times B\times \R^\krest}_{\R^\times }&\to
\bigl(t\cdot\frr\bigr)\dow^{A\times B\times \R^\krest}_{\R^\times };
\\
\frr_j\dow^{A\times Y\times \R^\krest}_{\R^\times }+
\frr_j^-\dow^{A\times  \R^\krest}_{ \R^\times }
&\to
\frr\dow^{A\times Y\times \R^\krest}_{\R^\times }+
\frr^-\dow^{A\times  \R^\krest}_{ \R^\times };
\\
\bigl(t\cdot \frr_j\bigr)\dow^{X\times B\times \R^\krest}_{\R^\times }+
\bigl(t\cdot \frr_j^+\bigr)\dow^{B\times  \R^\krest}_{ \R^\times }&\to
\bigl(t\cdot \frr\bigr)\dow^{X\times B\times \R^\krest}_{\R^\times }+
\bigl(t\cdot \frr^+\bigr)\dow^{B\times  \R^\krest}_{ \R^\times }.
\end{align*}

Next, let $(M,\mu)$ be a space with a continuous $\sigma$-finite measure.
We embed the group $\Gms(M)$ to $\Pol_\infty(M,M)$
assigning to $q\in\Gms(M)$ the pushforward $\frr_q$ of $\mu$ under the map $m\mapsto (m,q(m),q'(m))$
and assuming $\frr^+=0$, $\frr^-=0$.



\begin{conjecture}
	\label{con:1}
{\rm a)} The multiplication in $\Gms(M)$ extends by separate continuity to a multiplication%
\footnote{Cf. \cite{Ner-bist}, Subsect.3.7.}
in $\Pol_\infty(M,M)$. The latter semigroup extends to a category, whose objects are $\sigma$-finite 
Lebesgue spaces and sets of morphisms are $\Pol_\infty(X,Y)$.

\sm 
 
 {\rm b)} There is a counterpart of Theorem {\rm \ref{th:1}} for actions of groups on the space with infinite continuous measure.
 \end{conjecture}

3) According Nessonov \cite{Ness}, any non-singular action of $\ov S_\infty$ is equivalent to a measure preserving action
(on space with finite of $\sigma$-finite measures). So, for 
$K=\ov S_\infty$ we can   omit the condition '$K$ acts by measure preserving transformations' in Theorem \ref{th:1}. Apparently, this case is not unique. 

\begin{conjecture}
	\label{con:2}
	{\rm a)} A similar statement holds for complete  groups%
	\footnote{All these groups are 'heavy groups' in the sense of \cite{Ner-book}, Chapter VIII.} 
	$\ov {\mathrm{O}}(\infty)$, $\ov \U(\infty)$, $\ov \Sp(\infty)$ {\rm(}i.e., groups of all unitary operators in real, complex, quaternionic
	$\ell^2${\rm)}, for groups $\Ams[0,1]$, $\Ams_\infty(\R)$. 
	
	\sm 
	
	{\rm b)} A similar statement holds for all oligomorphic groups and inverse limits of oligomorphic groups.
\end{conjecture}

The most strong variant of this conjecture is: 

\sm

--- {\it the statement holds for all Roelcke precompact Polish groups.}

\sm 

On the definition of Roelcke precompact groups (there are 3 essentially different uniform structures on a topological group,
introduced by A.~Weil, D.~A.~Raikov, W.~Roelcke), see e.g., \cite{BIT}, \cite{BYa}; oligomorphic groups and groups
mentioned in claim a) are Roelcke precompact; it is strange enough that the group $\Gms(M)$ also is Roelcke precompact, 
see Ben Yaacov \cite{BYa}, Theorem 2.4.

 Math. Dept., University of Vienna; \\
  Institute for Information Transmission Problems;\\
 MechMath Dept., Moscow State University;\\
 yurii.neretin@univie.ac.at;\\
 URL: http://mat.univie.ac.at/$\sim$neretin/

\end{document}